\documentclass{amsart}
\usepackage{amssymb,mathrsfs}
\DeclareMathSymbol{\varGamma}{\mathord}{letters}{"00}
\DeclareMathSymbol{\varPi}{\mathord}{letters}{"05}
\DeclareMathSymbol{\varLambda}{\mathord}{letters}{"03}
\newcommand{\Cs}{\mathscr{C}}
\newcommand{\Ps}{\mathscr{P}}
\newcommand{\rd}{\mathrm d}
\begin{document}

\title[Compatibility of conformal and projective structures]{A criterion for compatibility of conformal and projective structures}

\author{Vladimir S.  Matveev}
\address{Mathematics Institute, Fakult\"at f\"ur Mathematik und Informatik,
Friedrich-Schiller-Universit\"at Jena, 07737 Jena, Germany}
\email{vladimir.matveev@uni-jena.de}
\author{Andrzej Trautman}
\address{Institute of Theoretical Physics, Ho\.za 69, 00681 Warszawa, Poland}
\email{andrzej.trautman@fuw.edu.pl}

\begin{abstract}
In a space-time \(M\), a conformal structure is defined by the distribution of light-cones.
Geodesics are traced by freely falling  particles,   and the collection of all unparameterized geodesics determines the projective structure of \(M\). The article contains a formulation of the  necessary and sufficient conditions for these structures to be compatible, i.e. to come from a metric tensor which is then unique up to a constant factor. The theorem applies to all dimensions and signatures.
\end{abstract}

\maketitle

\section{Introduction and  remarks on the history of the problem}

Hermann Weyl, in his early papers on `infinitesimal geometry' \cite{Weyl18,Weyl21},  described the two structures that underlie
the geometry and physics of the four-dimensional space-time. The propagation of light determines light cones; the collection of all such cones gives a conformal structure \(\Cs\) of Lorentzian signature. Weyl pointed out that gravitation is described by a linear connection: particles, freely falling in a gravitational field, trace unparameterised geodesics that define symmetric linear connections, but only up to `projective transformations'. The collection of all such projectively related connections is  a projective structure \(\Ps\) on a manifold. (Precise definitions are given in the next section.) A Riemannian metric of Lorentzian signature uniquely determines both these structures; Weyl has shown that  two metrics \(g\) and \(g'\)  give the same two structures (conformal and projective) if and only if  \(g'={\rm const}\cdot g\);  see
\emph{Satz} 1 in \cite{Weyl21}. Weyl did not, however, consider the problem of whether a given pair of conformal and projective structures come from one metric tensor. Simple examples show that, in general,  they do not.

The problem raised by Weyl  has attracted, over the years, a considerable interest among physicists.  Ehlers, Pirani and Schild
wrote, on this subject, an influential paper that was recently reprinted as a `Golden Oldie' \cite{EPS72}. These authors argue in favour of founding the geometry of space-time on its conformal and
 projective structures rather than on the `chronometric' approach of J. L. Synge \cite{Synge60,Synge62}. They formulate a necessary condition that the pair \((\Cs, \Ps)\) must satisfy in order to result from one metric tensor. Namely, according to this EPS
  condition, as it will be called here, the  null geodesics of the conformal geometry should be also geodesics, or autoparallels, as defined by the projective structure. Ehlers, Pirani and Schild
 formulate further conditions that the structures \(\Cs\)
and \(\Ps\) should satisfy so as to come from a unique, up to a constant factor, metric tensor. However, they do not give sufficient conditions for this to be the case. More comments on that paper and further references can be found in \cite{TrautmanOnEPS}.

In this paper we present  a theorem giving the necessary and sufficient
 conditions for compatibility of  conformal and projective structures.
The theorem is algorithmic in the sense that, to determine compatibility of \(\Cs\) and \(\Ps\), it suffices to compute a few simple
 expressions formed from the components of \(g\in\Cs\) and \(\varGamma\in\Ps\). Our result is also effective: if these two structures are compatible, then
 a simple integration suffices to find   the corresponding  metric tensor.

\section{Definitions and the theorem}

We consider smooth -- of class \(C^\infty\) -- manifolds and maps. All geometric objects
on an \(n\)-dimensional manifold are referred to local coordinates \((x^i)\),\; \(i=1,\dots,n\).

A \emph{conformal structure} on a manifold \(M\) is an equivalence class \(\Cs\)  of
 metric tensors \(g\) with respect to the following equivalence relation
\[
g \sim  g'\;\Longleftrightarrow \mbox{there is a function \(\varphi\) on \(M\) such that \(g'=g \exp
{2\varphi}\).}
\]
If \(g\in\Cs\), then \(\Cs\) can be denoted by \([g]\).
No assumption is made on the signature of the metric tensors;
they can be properly Riemannian.

Two  symmetric linear connections
 \(\varGamma=(\varGamma^{i}_{jk})\)  and   \(\varGamma'=({\varGamma'}^{i}_{jk})\) are said to be \emph{projectively equivalent} if their geodesics differ only by parametrisation.  Projective equivalence  is clearly an equivalence relation on the set of all  symmetric linear connections on \(M\). An equivalence class \(\Ps\) with respect to this relation is  called a {\em projective structure};  it is denoted by \([\varGamma]\) if it contains \(\varGamma\).

  Projective equivalence    can be formulated as the condition
 \begin{center}
\(
   \varGamma \sim \varGamma'\in\Ps\;\
   \Longleftrightarrow\; \mbox{there is a 1-form \(\psi\) so that \(\varGamma^{\prime i}_{jk}=\varGamma^{i}_{jk}+\delta^i_j\psi_k+
   \delta^i_k\psi_j.\)}
   \)
\end{center}
In this form it appears in \cite{Weyl21}, but the essence of this
 result was given already by Tullio Levi-Civita in his very first publication \cite{Levi-Civita},
written at the age of 23.

Tracy Thomas \cite{Thomas25} observed that, given two symmetric linear connections, it is easy to check whether they are projectively equivalent by  computing the traceless quantity \(\varPi(\varGamma)\), which is nowadays called the  Thomas symbol,
 \begin{equation*}
   \varPi^i_{jk}(\varGamma)=\varGamma^{i}_{jk}-{\frac{1}{n+1}}\delta^i_j\varGamma^p_{pk}-
   {\frac{1}{n+1}}\delta^i_k\varGamma^p_{pj},\quad n=\dim M.
 \end{equation*}
Namely, two symmetric linear connections are projectively equivalent, if and only if, their  Thomas symbols coincide,
 \begin{equation} \label{e:GG}  \varPi(\varGamma)=\varPi(\varGamma')\;\Longleftrightarrow\;\mbox{\(\varGamma\) and \(\varGamma'\) are projectively equivalent.}
 \end{equation}

Let \(\digamma(g)\) be the Levi-Civita connection defined by  \(g\). In local coordinates,
 \begin{equation*}
   \digamma^i_{jk}(g)=\tfrac{1}{2} g^{ip}(\partial_kg_{pj}+\partial_jg_{pk}-\partial_pg_{jk}),
 \end{equation*}
so that
\begin{equation}\label{e:dig}
 \digamma^i_{jk}(g\exp{2\varphi})=\digamma^i_{jk}(g)+ \delta^i_j\partial_k\varphi+\delta^i_k\partial_j\varphi-g^{ip}g_{jk}\partial_p\varphi,
\end{equation}
where \(\partial_k\varphi={\partial\varphi}/{\partial x^k}\), etc. If \(u=(u^i)\) is a null vector, \(g_{ij}u^iu^j=0\), then
\[
\digamma^i_{jk}(g\exp{2\varphi})u^j u^k-\digamma^i_{jk}(g)u^j u^k\;\|\; u^i
\]
so that (unparameterised) null geodesics are well defined in conformal geometry in the sense the null geodesics of \(g\) are reparametrised null geodesics of \(g \exp(2\varphi)\).

\vspace{1ex}

\noindent{\bf Definition.} {\em The conformal and projective structures \(\Cs\) and \(\Ps\) are said to be \emph{compatible} if there is \(g\in\Cs\) such that \(\digamma(g)\in\Ps\).}

\vspace{1ex}

 Given  \(g\in\Cs\) and \(\varGamma\in\Ps\), from \eqref{e:GG} one obtains
\begin{equation}\label{e:defeq}
\mbox{\(\Cs\) and \(\Ps\) are compatible \(\Longleftrightarrow\)
 \( \exists\;\varphi\) such that}\quad  \varPi(\digamma(g\exp{2\varphi}))=\varPi(\varGamma).
\end{equation}
Since the difference of two connection coefficients is a tensor, so is
\begin{equation*}
  T^i_{jk}\stackrel{\rm {def}}{=}\varPi^i_{jk}(\digamma(g)-\varGamma).
\end{equation*}
The components of this tensor depend on  the components of the metric tensor and their first derivatives and on the components of the linear connection.
Substituting \eqref{e:dig} into \eqref{e:defeq}, one infers that compatibility of \(\Cs\) and \(\Ps\) is equivalent to the existence of \(\varphi\) such that
\begin{equation}\label{e:Meq}
 T^i_{jk}- g_{jk}g^{ip}\partial_p\varphi+{\frac{1}{n+1}}\delta^i_j
 \partial_k\varphi+{\frac{1}{n+1}}\delta^i_k\partial_j\varphi=0.
\end{equation}
Let
\begin{equation}\label{7}
  T^i={\frac{n+1}{(n+2)(n-1)}}g^{jk}T^i_{jk}\quad\mbox{and}\quad T_i=g_{ij}T^j.
\end{equation}
By contraction of \eqref{e:Meq}  with \(g^{jk}\) one   obtains
\begin{equation}\label{e:phi}
  \partial_i\varphi=T_i
\end{equation}
Substituting \(\partial_i\varphi\) determined by \eqref{e:phi} and \eqref{7}  into
\eqref{e:Meq}, one obtains the following
 condition on \(g\) and  \(\varGamma\):
\begin{equation}\label{e:A}
  T^i_{jk}- g_{jk}T^i+{\frac{1}{n+1}}
  \delta^i_jT_k+{\frac{1}{n+1}}\delta^i_kT_j=0.
\end{equation}
Since the second partial derivatives of  \(\varphi\) commute, from \eqref{e:phi} one obtains
\begin{equation}\label{e:B}
  \partial_jT_i-\partial_iT_j=0.
\end{equation}

\vspace{1ex}
\noindent{\bf Theorem.}  {\em The conditions \eqref{e:A} and \eqref{e:B} are necessary and sufficient for local
compatibility of the conformal and projective structures, defined on \(M\) by \(g\) and \(\varGamma\), respectively. If, moreover, the first cohomology group of \(M\) vanishes, then there is global compatibility.}

\vspace{1ex}

\noindent{\it Proof.} The conditions are necessary because they were derived under the assumption of compatibility. Condition (\ref{e:B}) implies the existence of a local -- in a neighbourhood of every  point --  solution \(\varphi\) of
\eqref{e:phi}. Replacing now \(T_i\) in \eqref{e:A} by \(\partial_i\varphi\) one obtains that condition \eqref{e:Meq}  holds. If the first cohomology group of \(M\) vanishes, then the closed form \(T_i\rd x^i\) is exact and thus \(\varphi\) is defined all over \(M\).
 {\hfill\(\Box\)}

\section{A simple application}

Using the theorem one can confirm the existence   of  pairs \((\Cs, \Ps)\)
  that are incompatible even though the EPS  condition  holds.
Indeed, let \(g\) be a Lorentzian metric on an \(n\)-dimensional manifold \(M\),\;
 \(n\geqslant 3\), and \(\Cs=[g]\).  Given a vector field \((S^i)\) on \(M\), one considers the projective
structure \(\Ps=[\varGamma]\) such that
\begin{equation}\label{e:Gamma}
  \varGamma^i_{jk}=\digamma^i_{jk}(g)-S^i g_{jk}.
\end{equation}
If \(u^i\) is a null vector, \(g_{ij}u^iu^j=0\), then \((\varGamma^i_{jk}-\digamma^i_{jk}(g))u^j u^k=0\)
so that a null geodesic with respect to \(\Cs\) is also a geodesic with respect to \(\Ps\)
and the EPS condition is satisfied.

Computing now \(T^i_{jk}\) for \(\varGamma\) given by \eqref{e:Gamma}, one obtains
\(T^i=S^i\) and that the algebraic condition \eqref{e:A} is satisfied. Therefore,
the pair \((\Cs,\Ps)\) now under consideration is compatible if, and only if, the form
\(g_{ij}S^j\rd x^i\) is closed. In other words, to obtain a manifold with a pair
\((\Cs,\Ps)\) that satisfies the EPS condition but is incompatible, it suffices to take
\(\Cs\) containing  a Lorentzian metric \(g\) and \(\Ps=[\varGamma]\) given by
\eqref{e:Gamma}, where \(S^i\) is vector field with a non-integrable distribution of subspaces orthogonal to it.

\section{Concluding remarks}   \label{conclusions}
The result presented in this paper, though technically very simple, completes
a line of research initiated by Weyl and continued by physicists. Many mathematical
objects consist of two -- or more -- structures on one set, connected by a notion of compatibility. Conformal and projective structures on manifolds have a clear origin in physics and, for   this reason, their compatibility has attracted interest of theoreticians.

The left hand sides of  \eqref{e:A} and \eqref{e:B} are  tensors of the type given by the position of  their indices. Moreover, they are determined by \(\Ps\) and \(\Cs\), but do not depend on the representatives of these equivalence classes. One can consider   these tensors as a measure of noncompatibility of the projective and conformal structures.  These tensors could be used in the construction and study  of those nonmetric relativistic theories of space-time that  use      the conformal and projective structures as the principal building blocks.

It is worth noting here that a conformal structure can be easily reconstructed from the knowledge of the distribution of light cones. Indeed, if \(v\in T_xM\) is a null vector,
then \(g_{ij} v^i v^j  = 0\) is a linear equation for the components of the metric tensor and, by taking \(n(n+1)/2-1\) generic null  vectors at a point, one  obtains a system of linear equations whose solution space is one dimensional and gives the conformal structure at that point. The somewhat subtler
procedure of reconstructing a symmetric linear connection from the  set of all unparametrized geodesics can also be reduced to solving a system of linear equations; see \cite[\S 2.1]{Matveev} for details.

Closely related to the question considered here -- but much more difficult -- is the \emph{Roger Liouvelle problem} initiated in \cite{RLiouville}: given a system of differential equations
\begin{equation}\label{e:Liouv}
  \ddot x^i=\varLambda^i(x,\dot x),\quad i=1,\dots,n,\quad \dot x^i=\rd x^i/\rd t,
\end{equation}
to find the conditions on the functions \(\varLambda\) so that the solutions of
\eqref{e:Liouv} represent geodesics of a Levi-Civita connection. Recently, the problem has been solved, in two dimensions, by Robert Bryant, Maciej Dunajski and Michael Eastwood
\cite{BryantEastwoodDunajski}.

It is also worth noting that there exist projective structures \(\Ps\) such that there is no metric \(g\)
satisfying \(\digamma(g)\in\Ps\). Indeed, by the results of  \cite{EastwoodMatveev},  the existence of such a metric
 is equivalent to the existence of covariant constant sections
of a non-trivial vector bundle with connection. For almost all projective structures such  parallel sections   do not exist. From general theory there follows the existence of  complete systems of {\em differential invariants}, i.e. invariant algebraic expressions in the components of \(\varGamma\)  and its derivatives that determine whether  there exists    a metric corresponding to the  projective structure \([\varGamma]\),  see e.g. \cite{Nurowski}.

\section*{Acknowledgments}

The authors, who have never met in person, thank Pawe{\l} Nurowski for the initiative and encouragement to write this paper.

\end{document}